\documentclass[12pt,leqno]{amsart}

\newtheorem{prop}{Proposition}[section]
\newtheorem{thm}[prop]{Theorem}
\newtheorem{lem}[prop]{Lemma}
\newtheorem{cor}[prop]{Corollary}

\numberwithin{equation}{section}

\newcommand{\fh}{\mathfrak{h}}
\newcommand{\Q}{\mathbb{Q}}
\newcommand{\Z}{\mathbb{Z}}
\newcommand{\R}{\mathbb{R}}
\newcommand{\C}{\mathbb{C}}
\newcommand{\PP}{\mathbb{P}}
\newcommand{\co}{\mathcal{O}}

\DeclareMathSymbol{\ordcol}{\mathord}{operators}{'072}

\DeclareMathSymbol{\notdiv}{\mathbin}{AMSb}{"2D}
\newcommand{\hcs}{\raise0.37ex\hbox{\kern-1pt{\scshape\lowercase{c}}}}
\newcommand{\hc}{\raise0.53ex\hbox{\kern-1pt{{c}}}}

\begin{document}

\title{Principal moduli and Class Fields}

\author{David Cox} \address{Department of Mathematics and Computer
Science, Amherst College, Amherst, MA
01002-5000, USA}

\email{dac@cs.amherst.edu} 

\author{John M{\hcs}Kay}
\address{Department of Mathematics and Statistics, 1455 de
Maisonneuve West, Concordia University, Montr\'eal, Qu\'ebec, H3G 1M8,
Canada}

\email{mckay@cs.concordia.ca}

\author{Peter Stevenhagen} \address{Mathematisch Instituut,
Universiteit Leiden,
Postbus 9512,
2300 RA Leiden, The Netherlands}

\email{psh@math.leidenuniv.nl}


\begin{abstract} We study 
the values taken by $\Gamma_0(n)$-modular functions 
at elliptic points of order $2$ for the Fricke group
$\Gamma_0(n)^\dagger$ that lie outside $\Gamma_0(n)$.
In the case of a principal modulus (`Hauptmodul') for
$\Gamma_0(n)$ or $\Gamma_0(n)^\dagger$, we determine the
class fields generated by these values.
\end{abstract}

\maketitle

\section{Introduction}
\label{intro}

\noindent
The Fricke group $\Gamma_0(n)^\dagger$ of level $n > 1$ is
the subgroup of $\mathrm{PSL}_2(\R)$ generated by the
congruence subgroup $\Gamma_0(n)\subset \mathrm{PSL}_2(\Z)$
and the Fricke involution
$ 
w_n = \sqrt{n}
\left(\begin{smallmatrix}
0 & -\frac{1}{n}\\ 1 & \phantom{-}0\end{smallmatrix}\right),
$ 
which acts on the upper half plane $\fh$ by $w_n(\tau)=
\frac{-1}{n\tau}$. 
As explained in \cite[Section 4]{cn}, $\Gamma_0(n)^\dagger$ is a subgroup of
finite index of the full normalizer $\Gamma_0(n)^+$ of $\Gamma_0(n)$ in
$\mathrm{PSL}_2(\R)$, and equal to it if $n$ is prime.
The subgroup $\Gamma_0(n)$ has index~$2$ in $\Gamma_0(n)^\dagger$, and 
the elements of the coset $\Gamma_0(n)^\dagger\setminus \Gamma_0(n)$
are the matrices 
$\sqrt n
\left(\begin{smallmatrix} A & \frac{B}{n} \\ C & D\\\end{smallmatrix}\right)
(\mathrm{mod} \pm1)$
with $A, B, C, D\in\Z$ satisfying $nAD-BC=1$.
Such a matrix is elliptic of order 2 if and only if it has
trace $0$, i.e., if it is represented by
\begin{equation}
\label{ellelt}
\alpha = \sqrt{n}
\begin{pmatrix}
A & \!\!\phantom{-}\frac{B}{n}\\ C
& \!\!-A\end{pmatrix}\ \in \mathrm{SL}_2(\R)
\end{equation}
with $nA^2 + BC = -1$.
We always pick $\alpha$ in \eqref{ellelt} in such a way that 
we have $B < 0$ and $C>0$.
Then the fixed point of $\alpha$ in $\fh$ is the element
\begin{equation}
\label{tdef}
\tau_\alpha = \frac{nA+\sqrt{-n}}{nC} \in \fh,
\end{equation}
which is a zero of the polynomial
\begin{equation}
\label{fdef}
f_\alpha=nC X^2 -2nA X -B\in\Z [X]
\end{equation}
of discriminant $4n(nA^2+BC)=-4n$.
As $B$ is coprime to $nA$, the polynomial $f_\alpha$ is irreducible in
$\Z[X]$ unless $B$ and $C$ are both even.
In the latter case we have $n=(-1-BC)\frac{1}{A^2}\equiv 3\bmod 4$, and 
$\frac{1}{2}f_\alpha$ is irreducible in $\Z [X]$ of discriminant $-n$.
It follows (cf. \cite[Thm.\ 7.7 and Prop.\ 7.4]{cox}) that the complex lattice
\begin{equation}
\label{adef}
I_\alpha = [\tau_\alpha, 1]=\frac{1}{nC}[nA+\sqrt{-n},nC] \subset \C
\end{equation}
corresponding to the fixed point of $\alpha$ is an invertible
ideal for the quadratic order
$\co_\alpha=\{\lambda\in\C: \lambda I_\alpha\subset I_\alpha\}$,
which has discriminant $-n$ if $B$ and $C$ are both even,
and discriminant $-4n$ otherwise.
More explicitly, we have
\begin{equation}
\label{oa}
\co_\alpha = \begin{cases}
        \Z[\textstyle{\frac{-n+\sqrt{-n}}{2}}] &\text{for }n \equiv 3
        \bmod 4,\ B\ \text{and}\ C\ \text{even}\\ \Z[\sqrt{-n}] 
        & \text{otherwise.} \end{cases}
\end{equation}
The goal of this paper is to study the Galois theoretic properties of
the values $f(\tau_\alpha)$ over $\Q(\tau_\alpha)$ when $f$ is a
$\Gamma_0(n)$-modular function and $\alpha$ is an elliptic element of
order 2 in $\Gamma_0(n)^\dagger\setminus \Gamma_0(n)$.  If $f$ is
modular for the full modular group $\mathrm{PSL}_2(\Z)$ and has
rational $q$-expansion, i.e., if $f$ is an element of the field
$\Q(j)$ of modular functions of level 1, then the classical theory of
complex multiplication tells us that if $f(\tau_\alpha)$ is finite,
then it is contained in the ring class field of the order
$\co_\alpha$.  This class field is an abelian extension of
$\Q(\tau_\alpha)=\Q(\sqrt{-n})$, and its Galois group over
$\Q(\tau_\alpha)$ is canonically isomorphic under the Artin
isomorphism to the ideal class group $C(\co_\alpha)$ of the order
$\co_\alpha$.  We will write $\sigma_\alpha$ for the automorphism
corresponding to the ideal class $[I_\alpha]\in C(\co_\alpha)$.

For $f\in\Q(j)$ of level 1 as above, a classical result 
\cite[Cor. 11.37]{cox} on the Galois theoretic properties
of the values of $j$ at imaginary quadratic arguments
implies that we have
$\sigma_\beta\big(f(\tau_\alpha)\big) = f(\tau_{\gamma})$ when
$\alpha, \beta$ and $\gamma$ are elliptic elements as in \eqref{ellelt} 
satisfying $\co_\alpha = \co_\beta = \co_\gamma$ and
$[I_\alpha I_\beta^{-1}]=[I_\gamma]\in C(\co_\alpha)$.

Our main result below shows that $\Gamma_0(n)$-modular functions
behave like level 1 functions when evaluated at the elliptic points of
order 2 of the Fricke group $\Gamma_0(n)^\dagger$ that lie
outside $\Gamma_0(n)$.  The proof given in
Section~\ref{modfn} uses the Shimura reciprocity law
\cite[Thm.~6.31]{shimura} as explained in \cite{gs,stevenhagen}.

\begin{thm}
\label{mainthm}
Fix $n > 1$ and let $f$ be a $\Gamma_0(n)$-modular function with
rational $q$-expansion at $\infty$.  Also let $\alpha$ be an elliptic
element of order $2$ in $\Gamma_0(n)^\dagger\setminus \Gamma_0(n)$ as
in \eqref{ellelt}, and assume that $f$ is defined at $\tau_\alpha$.
Then:
\begin{enumerate}
\item $f(\tau_\alpha)$ lies in the ring class field of $\co_\alpha$. 
\item Let $\beta, \gamma$ be elliptic elements as in \eqref{ellelt} with
$\co_\alpha = \co_\beta = \co_\gamma$, and suppose
$[I_\alpha I_\beta^{-1}]=[I_\gamma]\in C(\co_\alpha)$.
Then
we have
\[
\sigma_\beta\big(f(\tau_\alpha)\big) = f(\tau_{\gamma}).
\]
\end{enumerate}
\end{thm}
\noindent
Our motivation for proving Theorem~\ref{mainthm} goes back to the
special case where $f$ is a principal modulus (`Hauptmodul') of the
kind that arises in the moonshine conjectures for the monster group
\cite{cn} and its generalizations to replicable functions \cite{cun}.
These principal moduli are known to correspond to genus zero subgroups
of $\mathrm{PSL}_2(\R)$ that lie between $\Gamma_0(n)$ and its
normalizer $\Gamma_0(n)^+$ for some $n$.  In Section~\ref{haupt} we
apply Theorem~\ref{mainthm} to this situation, and explain how our
results relate to those of Chen and Yui \cite{cy}.  The main
corollary, which is proved in Section~\ref{haupt}, is the following.

\begin{cor}
\label{maincor}
Fix $n > 0$ and let $f$ be a principal modulus with rational
$q$-expansion at $\infty$ for either $\Gamma_0(n)$ or
$\Gamma_0(n)^\dagger$.  Suppose that $f$ is defined at the fixed point
$\tau_\alpha$ from \eqref{tdef}.  Then $\Q(\tau_\alpha,
f(\tau_\alpha))$ is the ring class field of the order $\co_\alpha$.
\end{cor}

\section{Ideals, Quadratic Forms, and Conjugacy Classes}
\label{elliptic}

\noindent
This section collects what we need about the ideals $I_\alpha$
arising in \eqref{adef} and their classes in $C(\co_\alpha)$.
While part of this material is well known (see \cite{helling,maclachlan}),
we include the details for the convenience of the reader.

In \eqref{fdef}, we denoted by $f_\alpha$ the quadratic polynomial
in $\Z[X]$ having the fixed point $\tau_\alpha$ of the
elliptic element $\alpha$ as one of its zeroes.
Its homogeneous form
\begin{equation}
\label{quadform}
F_\alpha(x,y) = nCx^2-2nAxy-By^2 \in\Z[x,y]
\end{equation}
is a binary quadratic form of discriminant $-4n$.
Since we take $B < 0$ and $C > 0$, the form $F_\alpha$ is positive definite.
Depending on the case we are in in \eqref{oa},
either $F_\alpha$ or $\frac{1}{2}F_\alpha$
is primitive of discriminant $\mathrm{disc}(\co_\alpha)\in\{-n, -4n\}$.
Under the standard bijection \cite[Thm.\ 7.7]{cox} between the set of
$\mathrm{SL}_2(\Z)$-orbits of primitive positive definite binary
quadratic forms of discriminant $\mathrm{disc}(\co_\alpha)$ and
the class group $C(\co_\alpha)$,
this primitive form corresponds to the ideal class
$[I_\alpha]\in C(\co_\alpha)$.

\begin{lem}
\label{conjeq}
Let $\alpha$ and $\alpha'$ be as in \eqref{ellelt}.  Then the
quadratic forms $F_\alpha$ and $F_{\alpha'}$ are
$\mathrm{SL}_2(\Z)$-equivalent if and only if $\alpha$ and $\alpha'$
are conjugate in $\Gamma_0(n)^\dagger$.  If $\alpha$ and $\alpha'$ are
conjugate in $\Gamma_0(n)^\dagger$, there exists $\delta\in
\Gamma_0(n)$ with $\delta^{-1} \alpha \delta = \alpha'$.
\end{lem}

\begin{proof}
Let $M_{\alpha}$ denote the matrix representing the form $\alpha$
from \eqref{quadform}:
\[
M_{\alpha} = 
\begin{pmatrix} \phantom{-}nC &  -nA \\ -nA & -B\ \end{pmatrix} 
= \sqrt{n} 
\begin{pmatrix} \phantom{-}0 & 1 \\ -1 & 0 \end{pmatrix} 
\alpha.
\]
Since every $\delta \in \mathrm{SL}_2(\Z)$ satisfies
\[
\begin{pmatrix} 0 & \!\!-1 \\ 1 & \!\!\phantom{-}0 \end{pmatrix} 
\delta^{t} 
\begin{pmatrix} \phantom{-}0 & 1 \\ -1 & 0 \end{pmatrix} = 
\delta^{-1},
\]
we have
\[
\delta^{t} M_{\alpha} \delta = M_{\alpha'} \Longleftrightarrow 
\delta^{-1} \alpha \delta = \alpha'.
\]
When $F_{\alpha}$ and $F_{\alpha'}$ are equivalent via $\delta\in
\mathrm{SL}_2(\Z)$, reduction modulo~$n$ gives
\[
\delta^{t} \begin{pmatrix} 0 & \phantom{-}0 \\ 0 &
-B\\ \end{pmatrix} \delta \equiv \begin{pmatrix} 0 & \phantom{-}0 \\ 0 &
-B'\\ \end{pmatrix} \bmod n.
\]
If $\delta = \left(\begin{smallmatrix} a & b \\ c &
d \end{smallmatrix}\right)$, this easily implies
\[
-Bc^2 \equiv -Bcd \equiv 0 \bmod n.
\]
Since $\gcd(B,n) = \gcd(c,d) = 1$ (because $\det(\alpha) =
\det(\delta) = 1$), we see that $c \equiv 0 \bmod n$, so that $\delta
\in \Gamma_{0}(n)$.  Then the above equivalence implies that $\alpha$
and~$\alpha'$ represent conjugate elements in $\Gamma_0(n)^\dagger$.
 
Conversely, suppose we have $\alpha' = \pm \delta^{-1} \alpha \delta$
(remember that conjugacy classes are computed in
$\Gamma_{0}(n)^\dagger \subset \mathrm{PSL}_2(\R)$).  If $\delta$
represents an element of $\Gamma_{0}(n)$, then the above equivalence
shows that $F_{\alpha}$ and $\pm F_{\alpha'}$ are properly equivalent,
and the sign must be $+$ since $F_{\alpha}$ and $F_{\alpha'}$ are
positive definite.  The final statement of the lemma says that this is
the only case we need to consider.  Indeed, if $\delta$ represents an
element of $\Gamma_{0}(n)^\dagger \setminus \Gamma_{0}(n)$, then
$\alpha^{-1}\delta$ represents an element of $\Gamma_{0}(n)$ since
$\Gamma_{0}(n)$ has index $2$ in $\Gamma_{0}(n)^\dagger$, and we may
replace $\delta$ by $\alpha^{-1}\delta$:
\[
\alpha' = \pm \delta^{-1}\alpha\delta = \pm (\alpha^{-1}\delta)^{-1}
\alpha (\alpha^{-1}\delta).
\]
This completes the proof of the lemma.
\end{proof}
\noindent
It follows from Lemma \ref{conjeq} that $\Gamma_{0}(n)$-modular
functions have a well-defined value on the
$\Gamma_{0}(n)^\dagger$-orbits of the fixed points in \eqref{tdef},
and that elliptic elements $\alpha, \beta$ in \eqref{ellelt} with
$\co_\alpha=\co_\beta$ are conjugate in $\Gamma_{0}(n)^\dagger$ if and
only if $\sigma_\alpha$ and $\sigma_\beta$ coincide on the ring class
field of $\co_\alpha$.

Our main result (Theorem \ref{mainthm}) gives a Galois action in terms
of Artin automorphisms $\sigma_\alpha$ corresponding to the ideal
classes $[I_\alpha]$.  It is therefore useful to know that every ideal
class in $C(\co_\alpha)$ is of this form.

\begin{lem}
\label{giveall}
Let $n$ be a positive integer and fix $\co$ to be either
$\Z[\sqrt{-n}]$ or $\Z[\frac{n+\sqrt{-n}}{2}]$, where in the latter
case we assume $n \equiv 3 \bmod 4$.  Then the ideals $I_\alpha =
\frac{1}{nC}[nA+\sqrt{-n},nC]$ from\eqref{adef} with
$\co_\alpha = \co$ represent all ideal
classes in the ideal class group $C(\co)$.
\end{lem}

\begin{proof}
It suffices to show that a primitive positive definite quadratic form
$G$ of discriminant $D = -4n$, or $D=-n$ with $n \equiv 3 \bmod 4$, is
properly equivalent to $F_\alpha$ or $\frac12 F_\alpha$ respectively.
By \cite[Lemmas 2.3 and 2.25]{cox}, we may assume that $G(x,y) = ax^2
+ bxy + cy^2$ has $\mathrm{gcd}(c,D) = 1$.

First suppose that $D = -4n$.  Since $b$ is even in this case, we can
find an integer $k$ such that
\[
ck \equiv -b/2 \bmod n.
\]
Using this $k$ to replace $G(x,y)$ with the equivalent form 
\begin{equation}
\label{greduce}
G(x,y+kx) = (a+bk+ck^2)x^2 + (b+2kc)xy + cy^2,
\end{equation}
we may assume that $b$ is divisible by $2n$.  Then $-4n = b^2 - 4ac$
and $\mathrm{gcd}(c,-4n) = 1$ imply that $a$ is divisible by $n$.  It
follows that $G$ is of the form $F_\alpha$, as desired.

Next assume $D = -n$ and $n \equiv 3 \bmod 4$.  Then we can find
$k$ such that
\[
2ck \equiv -b \bmod n.
\] 
Using this $k$ and \eqref{greduce}, we may assume that $b$ is
divisible by $n$.  Since $n$ is odd, $-n = b^2 - 4ac$ and
$\mathrm{gcd}(c,-n) = 1$ imply that $a$ is divisible by $n$.  Then $G
= \frac12 F_\alpha$ follows easily.
\end{proof}
\noindent
For $\co=\Z[\sqrt{-n}]$, the form $F_\alpha=nx^2+y^2$
corresponding to the principal ideal generated by
$\tau=\frac{-1}{\sqrt{-n}}$
is in the unit class in $C(\co)$.
For $\co=\Z[\frac{-n+\sqrt{-n}}{2}]$, we can start from 
$G(x,y)=\frac{n+1}{4}x^2+xy+y^2$ in the unit class and apply the proof of
\eqref{giveall} (with $k=\frac{n-1}{2}$) to obtain the element
$F_\alpha(x,y)=\frac{n(n+1)}{4}x^2+nxy+y^2$ in the unit class.
It corresponds to the principal ideal generated by
$\tau=\frac{-2n+2\sqrt{-n}}{n(n+1)}$.
Note that in both cases, $\tau\co$ is the lattice $[\tau, 1]$.

\section{Proof of the Main Theorem}
\label{modfn}
\noindent
We now prove Theorem~\ref{mainthm} using the method explained
in \cite{stevenhagen}.

Let $f$ be a modular function for
$\Gamma_0(n)$ with rational $q$-expansion and assume that $f$ is
defined at the fixed point $\tau_\alpha$ from \eqref{tdef}. 
Our first task is to show
that $f(\tau_\alpha)$ lies in the ring class field of $\co_\alpha$. 
We first do this when the corresponding ideal $I_\alpha$ is in the
unit class in $C(\co_\alpha)$.
This leads to the following special choice for $\tau=\tau_\alpha$
and in each of the two cases for $\co=\co_\alpha$ provided by \eqref{oa}.
\[
\begin{tabular}{r | c | c | c | c}
& $n$ & $\alpha$ & $\tau$ & $\co$ \\ \hline
Case 1 & arbitrary & $\sqrt{n}\begin{pmatrix} 0 & \!\!-\frac1n\\ 1 &
\!\!\phantom{-} 
0\end{pmatrix}$ & $\frac{-1}{\sqrt{-n}}$ &
$\Z[n\tau]$ \\ \hline Case 2 & $n \equiv 3 \bmod 4$ &
$\sqrt{n}\begin{pmatrix} 1 & 
\!\!-\frac2n\\ \frac{n+1}{2} & \!\!-1\end{pmatrix}$ &
$\frac{-2n+2\sqrt{-n}}{n(n+1)}$ & $\Z[\frac{n(n+1)}{4}\tau]$
\end{tabular}
\]
Standard results in complex multiplication imply 
that $f(\tau)$ lies in the maximal Abelian extension $K^\mathrm{ab}$
of $K = \Q(\tau)= \Q(\sqrt{-n})$.
Class field theory describes the absolute abelian Galois group 
$\mathrm{Gal}(K^\mathrm{ab}/K)$ as the surjective image under the Artin
map of the group $\widehat{K}^*=(K\otimes_\Z \widehat{\Z})^*$ of
finite $K$-ideles.
(We generally write $\widehat A$ to denote the tensor product over
$\Z$ of a ring $A$ with the profinite completion
$\widehat{\Z}=\prod_p\widehat{\Z}_p$ of $\Z$.)
In order to show that $f(\tau)$ lies in the ring class field $H_\co$ of $\co$,
it suffices to show that every automorphism of $K^\mathrm{ab}$ over
$H_\co$ leaves $f(\tau)$ invariant.
As $\mathrm{Gal}(K^\mathrm{ab}/H_\co)$ is the image under the Artin map
of the subgroup $\widehat{\co}^*=(\co\otimes_\Z
\widehat{\Z})^*\subset \widehat{K}^*$, we need to show that
the Artin symbol of any idele $x\in \widehat{\co}^*$ leaves $f(\tau)$
invariant.

Shimura's reciprocity law tells us how (the Artin symbol of)
the idele $x\in \widehat{\co}^*$ acts on $f(\tau)$: we have
\begin{equation}
\label{srl}
f(\tau)^x= f^{g_\tau(x^{-1})}(\tau),
\end{equation}
where $g_{\tau} : \widehat{\co}^* \to \mathrm{GL}_2(\widehat{\Z})$
is the map that sends $x \in \widehat{\co}^*$ to the transpose of the
matrix of
multiplication by $x$ on the free $\widehat{\Z}$-module
$\widehat{\Z}\tau + \widehat{\Z}$ with respect to the basis
$\{\tau,1\}$.
One often writes 
$g_\tau(x)\big(\hskip-.7pt\genfrac{}{}{0pt}{}{\tau}{1}\hskip-.7pt\big) = 
\big(\hskip-1pt\genfrac{}{}{0pt}{}{x\tau}{x}\hskip-1pt\big)$. 
Note that $\mathrm{GL}_2(\widehat{\Z})$ acts in a natural
way on the field of modular functions over $\Q$ (see \cite[\S2 of Ch.\
  7]{lang}). 
For modular functions of level $N$, the action factors via
the quotient ${GL}_2({\Z}/N{\Z})$ (see \cite[\S3 of Ch.\
  6]{lang}).

{}From the irreducible polynomial of $\tau$ in $\Z[X]$ one easily
computes $g_\tau$ as in \cite[(3.3)]{stevenhagen} to be
\begin{equation}
\label{gtdef}
g_\tau(x) = g_\tau(a + bn\tau) = \begin{pmatrix} a & \!\!-b \\ nb
&
\!\!\phantom{-}a \end{pmatrix}
\end{equation}
for the order $\co=\Z[n\tau]$ in Case 1, and
\begin{equation}
\label{gtdef2}
g_\tau\big(a + b{\textstyle{\frac{n(n+1)}{4}}\tau}\big) =
\begin{pmatrix} a - nb& \!\!-b \\ \frac{n(n+1)}{4}b & 
\!\!\phantom{-}a \end{pmatrix} 
\end{equation}
for the order $\co=\Z[\frac{n(n+1)}{4}\tau]$ in Case 2.
In both cases, the matrix $g_\tau(x)$ is upper triangular modulo $n$
for every choice of $x\in \widehat{\co}^*$.
This means that modulo $n$, it is a product
\[
M_1 M_2 = \begin{pmatrix}1&0\\ 0&d\end{pmatrix} M_2 \in
\mathrm{GL}_2({\Z}/n{\Z}), 
\]
where $M_2$ is the reduction modulo $n$ of a matrix from
$\Gamma_0(n)$.  But $M_2$ acts trivially on $f$ since $f$ is a modular
function for $\Gamma_0(n)$, and $M_1$, which acts on $f$ via its Fourier
coefficients, also acts trivially since the $q$-expansion is rational. 
It follows that $f$ is invariant under
$\widehat{\co}^*$, so that $f(\tau)$ lies in the desired ring class field
$H_\co$.  
In particular, the class group $C(\co)$, which is naturally isomorphic to
the Galois group $\mathrm{Gal}(H_\co/K)$ under the Artin map, 
now acts on $f(\tau)$.

Our next task is to show that for elliptic elements $\beta$ and $\gamma$
satisfying $\co_\beta = \co_\gamma=\co$
and $[I_\beta^{-1}]=[I_\gamma]\in C(\co)$, the Galois action
of $\sigma_\beta$ on $f(\tau)$ is in each of the two cases given by
\begin{equation}
\label{action}
\sigma_\beta(f(\tau))=f(\tau_\gamma).
\end{equation}
This shows that we have $\sigma_\gamma(f(\tau_\gamma))=f(\tau)$ for
all $f(\tau_\gamma)$ with $\gamma$ as in \eqref{ellelt} having
$\co_\gamma=\co$.  In particular, all these values $f(\tau_\gamma)$
are conjugate over $K$ as soon as one of them is known to be finite,
and one easily derives the second statement in \eqref{mainthm} from
\eqref{action}.  We may therefore finish the proof of \eqref{mainthm}
by proving \eqref{action}.

To prove \eqref{action},
we apply \eqref{srl} once more, with $x\in\widehat{K}^*$
an $\widehat{\co}$-generator of
$\widehat{I_\gamma}=I_\gamma\otimes_\Z \widehat{\Z}$ and
$g_\tau: \widehat{K}^* \to \mathrm{GL}_2(\widehat{\Q})$ 
the natural $\Q$-linear extension of the map $g_\tau: \widehat{\co}^* \to
\mathrm{GL}_2(\widehat{\Z})$ we had before.
As in \cite[Section~6]{stevenhagen}, we let $M\in
\mathrm{GL}^+_2({\Q})$ be the transpose of a matrix mapping
the lattice $[\tau,1]=\tau\co$ to the lattice $\tau I_\gamma$.
Write $I_\gamma=\frac{1}{nC}[nA+\sqrt{-n},nC]$ as in \eqref{adef}.
In Case 1 we have $\tau\sqrt{-n}=-1$, so
\[
\tau I_\gamma=\tau[\textstyle{\frac{nA+\sqrt{-n}}{nC}},1]
             =[\textstyle{\frac{nA\tau-1}{nC}},\tau],
\]
and we can take $M\in \mathrm{GL}^+_2({\Q})$ to be
\[
M=
\begin{pmatrix}A/C&\!\!-1/nC\\ 1&\!\!\phantom{-}0\end{pmatrix} =
\begin{pmatrix}1/C&A/C\\ 0&1\end{pmatrix} 
\begin{pmatrix}0&-1/n\\ 1&0\end{pmatrix} .
\]
In Case 2 our element $\tau$ satisfies $\tau(n+\sqrt{-n})=-2$, 
so we have
\[
\tau I_\gamma=
\tau[\textstyle{\frac{nA+\sqrt{-n}}{nC}},1]=
[\textstyle{\frac{n(A-1)\tau-2}{nC}},\tau]
\]
and we can take the matrix $M\in \mathrm{GL}^+_2({\Q})$ to be
\[
M=\begin{pmatrix}(A-1)/C&-2/nC\\ 1&\!\!\phantom{-}0\end{pmatrix} =
\begin{pmatrix}2/C&(A-1)/C\\ 0&1\end{pmatrix}
\begin{pmatrix}0&-1/n\\ 1&0\end{pmatrix} .
\]
By construction, we have $M(\tau)=\tau_\gamma$ in both cases.
Note that at this point we use the fact that the lattice $[\tau, 1]$
is a principal $\co$-ideal.
Further $M$ and $g_\tau(x)$ are transposes of matrices mapping the
$\widehat{\Z}$-module $\widehat{\Z}\tau+\widehat{\Z}\subset \widehat K$
onto $\widehat{I_\gamma}$,
so the matrix $g_\tau(x)M^{-1}$ is in $\mathrm{GL}_2(\widehat{\Z})$
as its transpose fixes $\widehat{\Z}\tau+\widehat{\Z}$. 
We can now compute
\begin{equation}
\label{shim}
\sigma_\beta(f(\tau))=f^{g_\tau(x)}(\tau)=
f^{g_\tau(x)M^{-1}}(M(\tau))=f^{g_\tau(x)M^{-1}}(\tau_\gamma)
\end{equation}
by evaluating $g_\tau(x)M^{-1}\in \mathrm{GL}_2(\widehat{\Z})$ modulo
the level $n$ of $f$.
All we need to know for this are the $p$-components $x_p$ of
a generator $x\in \widehat{K}^*$ of~$\widehat{I_\gamma}$
for the primes $p$ dividing $n$. 
Finding a generator $x_p$ of the $(\co\otimes_\Z {\Z}_p)$-ideal
$I_\gamma\otimes_\Z {\Z}_p$ is easy at primes $p|n$, since 
for these $p$ the identity $BC = -1-nA^2$ implies that 
$C$ is a unit in $\Z_p$.
This yields
\begin{equation}
\label{fazp}
\begin{aligned}
I_\gamma\otimes_\Z \Z_p &=
\textstyle{\frac{nA+\sqrt{-n}}{nC}}\Z_p +\Z_p=
(A-\textstyle{\frac{1}{\sqrt{-n}}})\Z_p +\Z_p\\
&=\textstyle{\frac{-1}{\sqrt{-n}}}\Z_p+\Z_p=
  \textstyle{\frac{-1}{\sqrt{-n}}}\Z_p[\sqrt{-n}].
\end{aligned}
\end{equation}
In Case 1 we have $\co\otimes_\Z {\Z}_p=\Z_p[\sqrt{-n}]$ and we
can choose $x_p=\frac{-1}{\sqrt{-n}}=\tau$ for all $p|n$.
Applying \eqref{gtdef} then yields
\[
g_\tau(x_p)=\begin{pmatrix}0&\!\!-1/n\\ 1& \!\!\phantom{-}
0\end{pmatrix}\in \mathrm{GL}_2({\Q}_p)
\] 
at all $p|n$, so $g_\tau(x)M^{-1}$ has
upper triangular reduction
$$
\begin{pmatrix}1/C&A/C\\ 0&1\end{pmatrix}^{-1}=
\begin{pmatrix}C&-A\\ 0&1\end{pmatrix}
$$
modulo~$n$.
It therefore leaves $f$ invariant, and \eqref{shim} reduces to 
\eqref{action}.

In Case 2 only odd primes can divide $n$, so again
we have $\co\otimes_\Z {\Z}_p=\Z_p[\sqrt{-n}]$ for $p$ dividing $n$,
and we may take
$x_p=\frac{n+1}{4}\tau=\frac{1+\sqrt{-n}}{2}\frac{-1}{\sqrt{-n}}$
since $\frac{1+\sqrt{-n}}{2}$ is a $p$-adic unit if $p$
divides $n$.
Applying \eqref{gtdef2} now yields
\[
g_\tau(x_p)=
\begin{pmatrix}-1&-1/n\\ \frac{n+1}{4}&0\end{pmatrix}=
\begin{pmatrix}1&-1\\ 0&\frac{n+1}{4}\end{pmatrix}
\begin{pmatrix}0&\!\!-1/n\\ 1& \!\!\phantom{-}0\end{pmatrix}
\in \mathrm{GL}_2({\Q}_p)
\]
at all $p|n$, so $g_\tau(x)M^{-1}$ has
upper triangular reduction
\[
\begin{pmatrix}1&-1\\ 0&\frac{n+1}{4}\end{pmatrix}
\begin{pmatrix}2/C&(A-1)/C\\ 0&1\end{pmatrix}^{-1}
\]
modulo~$n$.
As in Case 1, we find that $g_\tau(x)M^{-1}$ leaves $f$ invariant.
This finishes the proof of Theorem \eqref{mainthm}.

By Lemma~\ref{giveall},  Theorem~\ref{mainthm} furnishes a
complete description of the Galois theoretic properties of 
$f(\tau_\alpha)$ over $\Q(\tau_\alpha)$.
This will be particularly useful in the cases in the next
section, where we know exactly
which elements in $\mathrm{PSL}_2(\R)$ fix $f$.

\section{Principal moduli}
\label{haupt}

\noindent
A modular function $f$ for a discrete group $G \subset
\mathrm{PSL}_2(\R)$ commensurable with $\mathrm{PSL}_2(\Z)$ is a
\emph{principal modulus} (in German:\ \emph{Hauptmodul}) if it
generates the field of all modular functions for $G$.  This implies
$G\backslash\fh^* \simeq \PP^1$ via $f$.  In particular, $G$ has genus
$0$.

Corollary \ref{maincor}, which we will prove now, states that if the
modular function $f$ in Theorem~\ref{mainthm} is a principal modulus
for $\Gamma_0(n)$ or $\Gamma_0(n)^\dagger$, then we get a primitive
element of the ring class field.

\begin{proof}[Proof of Corollary \ref{maincor}] The case $n = 1$ is
obvious, so we may assume $n > 1$.  Lemma~\ref{giveall} shows that the
ideals $I_\alpha$ represent all ideal classes in $C(\co_\alpha)$.
Then Lemma~\ref{conjeq} and the correspondence between ideals and
quadratic forms imply that we get $|C(\co_\alpha)|$ different
$\Gamma_0(n)^\dagger$-inequivalent points of the form $\tau_\alpha$.
If $f$ is a principal modulus for $\Gamma_0(n)^\dagger$, it follows
that the complex numbers $f(\tau_\alpha)$ have $|C(\co_\alpha)|$
distinct values as we vary $\alpha$, and the same is true if $f$ is a
principal modulus for $\Gamma_0(n)$ since
$\Gamma_0(n)^\dagger$-inequivalent points are automatically
$\Gamma_0(n)$-inequivalent.  Thus the minimal polynomial of
$f(\tau_\alpha)$ over $\Q(\tau_\alpha)$ has degree $|C(\co_\alpha)|$
by Theorem~\ref{mainthm}.  The theorem follows since this is the
degree of the ring class field over $\Q(\tau_\alpha)$.
\end{proof}
\noindent
It is well known that $\Gamma_0(n)$ has genus $0$ for the $15$ numbers
\begin{equation}
\label{gzn0}
n = 1\text{--}10, 12, 13, 16, 18, 25.
\end{equation}
Turning to $\Gamma_0(n)^\dagger$, Ogg \cite{ogg} showed that
$\Gamma_0(n)^\dagger$ has genus $0$ for the following $37$ values of
$n > 1$ (we exclude $n = 1$ since $\Gamma_0(1)^\dagger =
\Gamma_0(1)$):
\begin{equation}
\label{ntable}
n = 2\text{--}21, 
     23\text{--}27, 29, 31, 32, 35, 36, 
     39, 41, 47, 49, 50, 59, 71. 
\end{equation}
However, in order to use Corollary \ref{maincor} for either
$\Gamma_0(n)$ or $\Gamma_0(n)^\dagger$, we also need to know that the
group has a principal modulus with a rational $q$-expansion.  As we
are dealing with a finite list of groups, one can simply check the
tables in \cite{cn}, which give explicit formulas for principal moduli
in all cases when $\Gamma_0(n)$ or $\Gamma_0(n)^\dagger$ has genus
$0$.  These all have rational $q$-expansions.  Alternatively, as S.\
Norton pointed out to us, Theorem~1 of \cite{norton} easily implies
the existence of a principal modulus with rational $q$-expansion
whenever $\Gamma_0(n)$ or $\Gamma_0(n)^\dagger$ has genus $0$.  It
follows that Corollary \ref{maincor} applies to all values of $n$ in
\eqref{gzn0} for $\Gamma_0(n)$ and \eqref{ntable} for
$\Gamma_0(n)^\dagger$.

The primitive elements for the ring class fields obtained from
principal moduli are often much smaller than those provided by the
$j$-function.  For example, take $f$ to be the function listed as 71A
in \cite[p.~337]{cn}, which is a principal modulus for the Fricke
group $\Gamma_0(71)^\dagger$.  There are 14 elliptic fixed points of
order 2 for $\Gamma_0(71)^\dagger$, of which seven are of discriminant
$-71$ and seven are of discriminant $-4\cdot71=-284$.  This is in
accordance with the fact that the two quadratic orders having these
discriminants have class number 7.  In both cases the corresponding
ring class field is the Hilbert class field $H$ of $\Q(\sqrt{-71})$.

When $n = 71$, the two special elliptic points of order 2 used in the
proof of Theorem~\ref{mainthm} are
\[
\tau_{-284} = \frac{-1}{\sqrt{-71}}\quad \text{and}\quad
\tau_{-71} = \frac{-71 + \sqrt{-71}}{71\cdot 36}.
\]
corresponding to the orders of discriminant $-284$ and $-71$
respectively.  One can compute that the minimal polynomial of
$f(\tau_{-284})$ is 
\[
h_{-284} = x^7-7x^5-11x^4+5x^3+18x^2+4x-11,
\]
while the minimal polynomial of $f(\tau_{-71})$ is
\[
h_{-71} = x^7+4x^6+5x^5+x^4-3x^3-2x^2+1.
\]
Each of these polynomials is a small generating polynomial for $H$
over $\Q(\sqrt{-71})$ and compares favorably with the polynomial
\[
w_{-71} = x^7-2x^6-x^5+x^4+x^3+x^2-x-1
\]
found by Weber \cite[Vol.\ III, p.\ 723]{weber} (this is the minimal
polynomial of $\mathfrak{f}(\sqrt{-71})/\sqrt{2}$, where
$\mathfrak{f}$ is the Weber function defined on page 114 of
\cite[Vol.\ III]{weber}).  If $\beta$ is a root of $w_{-71}$, then an
easy calculation shows that
\begin{align*}
\beta^2 -1-\beta^{-1}\ \ &\text{is a root of}\
 h_{-284}\\
-\beta^6+3\beta^5 -2\beta^4+ 1\ \ &\text{is a root of}\
 h_{-71}.
\end{align*}
Note that the values $f(\tau_{-71})$ and $f(\tau_{-284})$ are in this
case \emph{integral}, a general phenomenon for which we currently do
not have a proof.  For the importance of $h_{-284}$ and $h_{-71}$ for
the Schwarzian differential equation satisfied by~$f$, see \cite{hm}.

As we mentioned in the introduction, the moonshine conjectures for the
monster group provide us with principal moduli for genus zero
subgroups of $\mathrm{PSL}_2(\R)$ in the form of so-called
M{\hc}Kay-Thompson series.  Their singular values have been studied by
Chen and Yui in \cite{cy}.  They restrict to fundamental
M{\hc}Kay-Thompson series, i.e., the series $T_g$ associated to a
conjugacy class $g$ of the monster group that is the principal modulus
for a subgroup of $\mathrm{PSL}_2(\R)$ whose level $n$ equals the
order of~$g$.  If $T_g$ is fundamental, $\tau \in \fh$ is a CM-point
and $\co=\Z[a\tau]$ is the multiplier ring of $[\tau, 1]$, they prove
the following in \cite[Theorem~3.7.5]{cy}:
\begin{enumerate}
\item If $T_g$ is a principal modulus for $\Gamma_0(n)$,
then $\Q(\tau, T_g(\tau))$ is the ring class field of the order of index
$\frac{n}{\mathrm{gcd}(a,n)}$ in $\co$.
\item If $T_g$ is a principal modulus for $\Gamma_0(n)^\dagger$ with $n$ prime
and coprime to $a$,
then $\Q(\tau, T_g(\tau))$ is the ring class field of the order of index~$n$
in $\co$.
\end{enumerate}
Applying the first part of this theorem to $T_g(\tau_\alpha)$, one can
prove Corollary~\ref{maincor} when $T_g$ is a principal modulus for
$\Gamma_0(n)$.  So for $T_g(\tau_\alpha)$, the Chen-Yui result covers
all numbers in \eqref{gzn0}.  On the other hand, if $T_g$ is a
principal modulus for $\Gamma_0(n)^\dagger$, then the hypothesis
$\mathrm{gcd}(a,n) = 1$ means that the second part of Theorem~3.7.5
doesn't apply to $T_g(\tau_\alpha)$.  So for $T_g(\tau_\alpha)$, the
Chen-Yui result covers none of the numbers in \eqref{ntable}.  The
integrality results in \cite{cy} for the numbers $T_g(\tau)$, although
interesting in their own right, do not apply to the fixed points
$\tau=\tau_\alpha$ in this paper.  Numerical computations do however
suggest that $T_g(\tau_\alpha)$ is an algebraic integer, so this
question deserves further study.

\section*{Acknowledgements}

\noindent
The research of the second author is partially supported by NSERC.
The first author is grateful to Heng Huat Chan for bringing
\cite{gs,stevenhagen} to his attention, and the second author would
like to thank Simon Norton for suggesting the relevance of
\cite{norton}.

\end{document}